\documentclass{amsart}
\newtheorem{theorem}{Theorem}[section]
\newtheorem{lemma}{Lemma}[section]

\theoremstyle{definition}

\theoremstyle{remark}

\numberwithin{equation}{section}

\begin{document}

\title[Growth estimates of meromorphic solutions of linear ODEs]
{Estimates on the Growth of Meromorphic Solutions of Linear Differential Equations with density conditions}

\author{Yik-Man Chiang}
\address{Department of Mathematics, Hong Kong University of Science and Technology, Clear Water Bay, Hong Kong, China}
\email{machiang@ust.hk}
\thanks{The research is partially supported by the University Grants Council of the Hong Kong Special Administrative Region, China (GRF601111)}



\date{December 20, 2013
}

\dedicatory{This paper is dedicated to the seventieth birthday of Ilpo Laine.}


\begin{abstract}
We give an alternative and simpler method for getting pointwise estimate of meromorphic solutions of homogeneous linear differential equations with coefficients meromorphic in a finite disk or in the open plane originally obtained by Hayman and the author. In particular, our estimates generally give better upper bounds for higher order derivatives of the meromorphic solutions under consideration, 
are valid, however, outside an exceptional set of finite logarithmic density.
 The estimates again show that the growth of meromorphic solutions with a positive deficiency at $\infty$ can be estimated in terms of initial conditions of the solution at or near the origin and the characteristic functions of the coefficients. 
\end{abstract}

\maketitle




\section{Introduction and main results}

We consider meromorphic solutions of the differential equation
\begin{equation}
	\label{E:ode}
		y^{(n)}(z) + \sum_{\nu=0}^{n-1} f_\nu (z)\, y^{(\nu)}(z)=0,
\end{equation}We apply freely the classical Nevanlinna Theory notation throughout this paper \cite{Ha75}, \cite{La93}.
where the coefficients $f_\nu (z)\, \ 0\le \nu\le n-1$ are meromorphic in $\mathbb{C}$. We mention that  Heittokangas, Korhonen and R\"atty\"a  obtained sharper estimates for analytic solutions when the coefficients are analytic functions in \cite{HKR2004}. They also considered non-homogeneous equations in \cite{HKR2009}. For an up-to-date account on the growth of meromorphic solutions of algebraic differential equations with meromorphic coefficients, we refer the reader to Hayman \cite{Ha1996}.
\medskip

Bank  asked,  if as in the case when the (\ref{E:ode}) admits an entire solution \cite{Ba75}, an meromorphic solution of (\ref{E:ode}) can be estimated in terms of growth of Nevanlinna characteristics of the meromorphic coefficients alone. In general, he showed that this statement is not true in  \cite{Ba76} by constructing an example for any given real-valued increasing function $\Phi(r)\uparrow +\infty$  on $(0,+\infty)$, then one can construct a first order linear differential equation with entire coefficient of zero Nevanlinna order such that the differential equation admits a meromorphic function $f$ with $T(r,\, f)>\Phi(r)$ as $r\to +\infty$. One would need some extra terms to bound the growth of meromorphic solutions. An example of such a result is given by Bank and Laine.

\begin{theorem}[\cite{BL77}]
Suppose that the coefficients of (\ref{E:ode}) are arbitrary
meromorphic functions and
that $y(z)$
is a meromorphic solution of (1.1).  If
\[
 \Phi(r)=\max_{0 \leq i \leq n} \Big(\log r, T(r,f_i)\Big),
\]
then for any $\sigma > 1$, there exist positive constants $c, c_1$ and
$r_0$, such that for $r \geq {r_0},$
\begin{equation}
	\label{E:BL-upper}
  T(r,y) \leq c \Big\{rN(\sigma r, y)+r^2 \exp 
  \Big(c_1 J(\sigma r) \log (rJ(\sigma r)\Big) \Big\},
\end{equation}
where
\[
  J(r)={\overline N} (r, 1/y) +\Phi (r).
\]
\end{theorem}
\medskip

We note that one needs counting function $\overline{N}(r,\, 1/y)$ of \textit{distinct zeros} in the $J(r)$ above as part of the upper bound in (\ref{E:BL-upper}).
\medskip

Since the equation (\ref{E:ode}) is linear, so one can deduce from the expression
	\[
		\frac{f^{(n)}(z)}{f(z)} + \sum_{\nu=1}^{n-1} f_\nu (z) \frac{f^{(\nu)}(z)}{f(z)}+f_0(z)=0,
	\]
that
	\[
		\overline{N}(r,\, f)\le \sum_{\nu=0}^{n-1} \overline{N}(r,\, f_\nu)
		\le \sum_{\nu=0}^{n-1} {T}(r,\, f_\nu),
	\]
indicating that the coefficients can only bound the \textit{distinct} poles of $f$. Indeed, the example constructed by Bank \cite{Ba76} mentioned above has poles of rapidly increasing multiplicities, that is $N(r,\, f)/\overline{N}(r,\, f)$ is unbounded. Hayman and the author \cite{ChiangHayman} showed that one can still bound the growth of a meromorphic solution $f$ of (\ref{E:ode}) in terms of the characteristic functions of coefficients alone if the solution $f$ has \textit{relatively few poles}. In particular, this means that $\delta(\infty,\, f)>0$.  This follows from the following result.
\medskip

\begin{theorem}[\cite{ChiangHayman}]\label{T:ChiangHayman}
Suppose that $0<\rho< r< R$ and suppose that the coefficients $f_\nu,\ 0\le\nu\le n-1$ of (\ref{E:ode}) are analytic on the path $\Gamma=\Gamma (\theta_0,\, \rho,\, t)$ defined by
 the segment
\[
\Gamma_1: z=\tau e^{i\theta_0},\quad \rho\le \tau\le t,
\]
followed by the circle
\[
\Gamma_2: z=te^{i\theta},\quad \theta_0\le\theta\le \theta_0+2\pi.
\]

We suppose that $y(z)$ is a solution of the equation (1.1) and define
\[
K=2\max\biggl\{1, \sup_{0\le\nu\le n}|y^{(\nu)}(z_0)|\biggr\},
\]
where $z_0=\rho e^{i\theta_0}$. We also define
\[
C=C(f_\nu,\,\rho ,\, r,\, R)= (n+2)\exp\biggl\{ 
\frac{20 R}{R-r}\sum_{\nu=0}^n T(R,\, f_\nu)+\biggl(\sum_{\nu=0}^n p_\nu\biggr)\log\biggl(\frac{R}{\rho}\biggr)
\biggr\},
\]
where $p_\nu$ is the multiplicity of the order of pole of $f_\nu$ at $z=0$. 
Then we have for $|z|=t$, where $t$ is some number such that $r<t<\frac{1}{4}(3r+R)$,
\[
|y^{(\nu)}(z)|< KC^\nu e^{(2\pi+1)CR},\quad 0\le \nu\le n.
\]
\end{theorem}
\bigskip

One can easily deduce when $R=+\infty$, and for a transcendental meromorphic $f$ with $\delta(\infty,\, f)>0$, then for $0<\varepsilon<\delta$, we have 
	\[
		T(r,\, y)\le \Big(\frac{1}{\delta-\varepsilon}\Big)\,(2\pi+1)\,R\, C.
	\]
\medskip

The main purpose of this paper is to give a shorter proof of a slightly different statement to Theorem \ref{T:ChiangHayman} and asymptotic results outside some exceptional sets using a different method. On the other hand, the original Theorem \ref{T:ChiangHayman} can deal with non-homogeneous (\ref{E:ode}), while our alternative can only deal with the (\ref{E:ode}). We prove
\medskip

\begin{theorem}\label{T:main}
Let $y$ be a meromorphic solution to the differential equation (\ref{E:ode}) with meromorphic coefficients $f_\nu,\, \nu=0, \cdots, n-1$ in $|z|=r< R\le +\infty$ such that $f_\nu$ has a pole of order $q_\nu\ge 0\ (0\le \nu\le n-1)$. Given a constant $C>1$ and $0<\eta< 3\,e/2$ and $r=|z|$ is outside a union of discs centred at the poles of $y$ such that the sum of radii is not greater than $4\eta R$, then there is a $B=B(C)>1$ and a path 
	\begin{equation}
		\Omega=\Omega(\theta,\, \rho,\, r)
	\end{equation}
consists of the line segment
\[
\Omega_1 : z =\tau e^{i\theta_0},\quad 0\le\rho\le \tau\le r
\]
followed by the circle
\[
\Omega_2 : z =re ^{i\theta},\quad \theta_0\le \theta< \theta_0+2\pi,
\]
on which the coefficients $f_\nu$ are analytic and we have, for $z$ on $\Omega$, 
\[
	\label{E:main-estimate}
\sum_{j=0}^{n-1}|y^{(j)}(z)|\le K_1\, e^{(2\pi+1)\,D\, r}\le K_1\, e^{(2\pi+1)\,D\,R},
\]
where
\[
	\label{E:initial}
K_1= \sum_{j=0}^{n-1}|y^{(j)}(z_0)|,\qquad z_0=\rho\, e^{i\theta_0}
\]
and
\begin{equation}
	\label{E:D}
		\begin{split}
D: & =D(f_\nu,\, \rho,\, r,\, R;\, \eta,\, B,\, C)\\
&\quad =n\bigg\{1+ \Big(\frac{R^{H(\eta )(\frac{R+2\,er}{R-2\,er})}}{r}\Big)^q\,\exp\Big[B\,(1+H(\eta))\Big(\frac{R+2\,er}{R-2\,er}\Big)\, T(CR)\Big]\bigg\}
		\end{split}
\end{equation}
and where 
	\[
\label{E:defn}
		T(r)=\max_{0\le\nu\le n-1} T(r,\, f_\nu),\quad q=\max_{0\le\nu\le n-1} {q_\nu}, \quad H(\eta)=2+\log\frac{3e}{2\eta}.
	\]
\end{theorem}
\bigskip

For any $r^\prime$, we choose $r$ outside a union of discs centred at the zeros of $y$ such that the sum of radii is not greater than $4\eta R$ such that $r^\prime <r <R$ as described in the Theorem \ref{T:main}, we have
	\[
		\begin{split}
			T(r^\prime,\, y) &\le N(r,\, y)+m(r,\, y)\\
			&\le 	 N(r,\, y)+ (2\pi+1)\, R\, D(f_\nu,\, \rho,\, r,\, R;\, \eta,\, B,\, C).
		\end{split}
	\]
\bigskip

This improves upon Bank and Laine's estimate mentioned above. Suppose that $\delta(\infty,\, y)>0$, we choose $r$ outside a union of discs centred at the zeros of $y$ such that the sum of radii is not greater than $4\eta R$ and sufficiently large such that $N(r,\, y)<(1-\delta+\varepsilon/2)\, T(r,\, y)$. Without loss of generality, we may also assume that $r$ is so chosen such that $|\log K_1|<\frac12\varepsilon T(r,\, y)$ so that $T(r,\, y) < (1-\delta+\varepsilon)\, T(r,\, y)+(2\pi+1)R\, D$. We can easily deduce 
	\[
		T(r^\prime,\, y)\le \Big(\frac{1}{\delta(\infty,\, y)-\varepsilon}\Big)\, (2\pi+1)R\, D.
	\]
\bigskip

\begin{theorem}\label{T:density}
Let $y(z)$ be a meromorphic solution to equation (\ref{E:ode}), and we choose $0< \eta <(1+\log 2)/(16\,e^{5/2})<1$. Then there is a constant $B>1$ such that for every $\varepsilon>0$ be given, there is a $r(\varepsilon)>0$  we have 
\medskip
	\[
		\label{E:density-bound}
			\log m(r,\, y^{(j)})\le 5B\big(1+H(\eta)\big)\, T(3\,e^2r)+[(5H(\eta)-1)q+1+\varepsilon]\, \log r
	\]
\medskip
$j=0,\, 1,\, \cdots, n-1$ holds for all $r>r(\varepsilon)$ except perhaps for a set of positive logarithmic density $\displaystyle\frac{16\eta e^{5/2}}{1+\log 2}$.
\end{theorem}
\bigskip

This result is to be compared with the following  density-type result also obtained previously in \cite{ChiangHayman}:

\begin{theorem}
	Let $y(z)$ be a meromorphic solution of (\ref{E:ode}) such that the  $f_\nu$ are not all constant, we have
		\begin{equation}
			\log m(r,\, y) < \Big(\sum_{\nu=0}^{n-1} T(r,\, f_\nu)\Big)\Bigg[(\log r) \log \Big(\sum_{\nu=0}^{n-1} T(r,\, f_\nu)\Big)\Bigg]^\sigma,
		\end{equation}
where $\sigma>1$ is a constant, to hold outside an exceptional set of finite logarithmic measure.
\end{theorem}
\bigskip

\section{Preliminaries} \label{S:pre}

Let us write ${\bf y}(z)=(y_0,\, \cdots,\, y_{n-1})^T$ where $y_j(z),\, j=0,\, \cdots, n-1$ are complex functions of $z$. We define $\|{\bf y}\|=\sum_{j=0}^{n-1}|y_j|$. Suppose further that ${\bf A}=(a_{ij}(z))$ is a square matrix then we define $\|{\bf A}\|=\sum_{i,j}|a_{ij}|$. We note that 
	\[
		\left\|\int {\bf A}\, dt\right\|\le \int \left\|{\bf A}\right\|\, |dt|
	\]
(see e.g., \cite[pp. 1--4]{Bellman}). 
\bigskip

\begin{lemma}\cite[pp. 21--22]{Levin}\label{L:lower-bound} Let $R>0$ and $f(z)$ be analytic in $|z|\le 2\,e\, R$ with $f(0)=1$, and let $\eta$ be an arbitrary positive constant not exceeding $3\,e/2$. Then we have
	\begin{equation}\label{E:lower-bound}
		\log|f(z)| > -H(\eta) \log M(2\,e\,R,\, f)
	\end{equation}
for all $z$ in $|z|\le R$ but outside a union of disks centred at the zeros of $f$ such that the sum of radii is not greater than $4\eta\, R$, where
	\[
		H(\eta)=2+ \log \frac{3\,e}{2\,\eta}.
	\]
\end{lemma}
\medskip

We also need the following quotient representation of meromorphic functions due to Miles \cite{Miles1972} and Rubel \cite[Chapter 14]{Rubel1996}.
\medskip

\begin{lemma}[\cite{Miles1972}]\label{L:Miles} Let $f$ be a meromorphic function in the plane, and let $C>1$ be a given constant, then there exist entire functions $f_1$ and $f_2$, and a constant $B=B(C)>0$ such that
\[
f(z)=\frac{f_1(z)}{f_2(z)},\quad\quad \text{and}\quad\quad T(r,\, f_j)\le B\,T(C\,r,\, f),
\]
$j=1,\, 2$ and $r>0$. Here both the constants $B$ and $C$ are absolute constants, i.e., they are independent of the function $f$.
\end{lemma}
\bigskip

\section{Proof of Theorem \ref{T:main}}

We state and prove our main lemma that leads to the proof of the Theorem \ref{T:main}. 
\begin{lemma}\label{L:on-path}  
Let $y$ be a meromorphic solution to the differential equation (\ref{E:ode}) with meromorphic coefficients $f_\nu,\, \nu=0, \cdots, n-1$ in $|z|=r< R\le +\infty$. Suppose that the coefficients $f_\nu,\, \nu=0, \cdots, n-1$ are analytic on the path $\Omega=\Omega(\theta,\, \rho,\, r)$ as defined in the Theorem \ref{T:main}. Suppose $z_0=\rho e^{i\theta_0}$, then for all $z$ on $\Omega$,
\begin{equation}
\sum_{j=0}^{n-1}|y^{(j)}(z)|\le K_1\, \exp\bigg[\Big(\max_{\Omega}\sum_{\nu=0}^{n-1} (|f_\nu(z)|+1)\Big)\,(2\pi+1)\,r\bigg],
\end{equation}
where $K_1$ is given in (\ref{E:initial}).
\end{lemma}

\bigskip
\begin{proof} 

It is well-known that equation (\ref{E:ode}) can be written in the matrix form
\begin{equation}
	\label{E:matrix-eqn}
{\bf F'}(z)={\bf A}(z)\, \textbf{F}(z),
\end{equation}
where ${\bf F}=(y,\, y',\, \cdots, y^{(n-1)})^T$, and
\begin{equation}
{\bf A}(z)=
\left(
\begin{matrix}
               0       & 1       & 0      & \cdots & 0 \\
               0       & 0       & 1      & \cdots & 0 \\
               \vdots  & \vdots  & \vdots & \ddots & \vdots\\
               0       & 0       & \cdots & \cdots & 1\\
               -f_0    & -f_1    & \cdots & \cdots & -f_{n-1}\\
\end{matrix}
\right).
\end{equation}
\bigskip

A solution to the above matrix equation (\ref{E:matrix-eqn}) is given by
\begin{equation}
{\bf F}(z)={\bf F}(z_0)+\int_{z_0}^z {\bf A}(t)\, \mathbf{F}(t)\, dt.
\end{equation}
\bigskip

We now apply Gronwall's inequality \cite[pp. 35--36]{Bellman} to (\ref{E:matrix-eqn}) to obtain
	\begin{equation}\label{E:matrix}
		\begin{split}
			\|{\bf F}(z)\| & \le  \|{\bf F}(z_0)\|+\int_{z_0}^z \|{\bf A}(t)\|\,\|{\bf                    F}(t)\|\,|dt|\\
               & \le \|{\bf F}(z_0)\|\,\exp\biggl(\int_{z_0}^z \|{\bf                    A}(t)		\|\,|dt|\biggr)\\
               & \le  \|{\bf F}(z_0)\|\, \exp\bigg[\Big(\max_{\Omega}\sum_{\nu=0}^{n-1} |f_\nu(z)|+(n-1)\Big)\,(2\pi+1)\,r\bigg]\\
               & <  \|{\bf F}(z_0)\|\, \exp\bigg[\Big(\max_{\Omega}\sum_{\nu=0}^{n-1} (|f_\nu(z)|+1)\Big)\,(2\pi+1)\,r\bigg],\\
	\end{split}
	\end{equation}
where we have parametrized the path $\Omega$ with respect to arc length. Clearly the length of $\Omega$ is $(2\pi+1)\,r$ at most. This proves Lemma \ref{L:on-path}.
\end{proof}
\bigskip

We are ready to prove the Theorem \ref{T:main}, which is a direct application of the Lemma \ref{L:on-path} and the two lemmas stated in \S\ref{S:pre} . 
\medskip

\subsubsection*{Proof of the Theorem \ref{T:main}}
Given $C>1$ be given. Then Miles' result in Lemma \ref{L:Miles} asserts that we can choose a $B>0$ such that we can write the coefficients in  $f_\nu=f_{\nu,1}/f_{\nu, 2}$ from (\ref{E:ode}) such that
	\[
		T(r,\, f_{\nu,\, j})\le B\, T(C\, r,\, f_\nu),\qquad  j=1,\, 2;\quad 0\le \nu\le n-1
	\]
for $r>0$. We first assume that $f_{\nu, 2}(0)\not=0$ for $ 0\le \nu\le n-1$, then it is easy to see that we may assume that $f_{\nu, 2}(0)=1$
after dividing the numerator and the denominator by a suitable constant. Lemmas \ref{L:lower-bound} and \ref{L:Miles} with 
	\begin{equation}\label{E:miles-1}
		\begin{split}
				\log|f_\nu| & =\log |f_{\nu,1}| + \log |f_{\nu, 2}|^{-1}\\
                 & \le \log M(r,\, f_{\nu,1})+ H(\eta)\log M(2\,e\,r,\, f_{\nu,2})\\
                  & \le \log^+ M(r,\, f_{\nu,1})+ H(\eta)\log^+ M(2\,e\,r,\, f_{\nu,2})\\
                 & \le \Big(\frac{R+r}{R-r}\Big)\, T(R,\, f_{\nu,1})+H(\eta)\,\Big(\frac{R+2\,e\,r}{R-2\,e\,r}\Big)\,                     T(R,\, f_{\nu,\, 2})\\
                 & \le B\, \Big(\frac{R+r}{R-r}\Big)\,T(CR,\,f_{\nu})+BH(\eta)\,\Big(\frac{R+2\,e\,r}{R-2\,e\,r}\Big)\,                T(CR,\, f_{\nu})\\
                 & \le B\big(1+H(\eta)\big) \Big(\frac{R+2\,e\,r}{R-2\,e\,r}\Big)\, T(CR, f_{\nu})\\
                 & \le B\big(1+H(\eta)\big)\, \Big(\frac{R+2\,e\,r}{R-2\,e\,r}\Big)\, T(CR),
		\end{split}
	\end{equation}
where
	\[
		T(r)=\max_{0\le \nu\le n-1} T(r,\, f_\nu).
	\]
\bigskip

If, however, $f_{\nu, 2}$ has a zero of order $q_\nu$ at $z=0$, we consider
\[
f_\nu=\frac {f_{\nu,1}/z^{q_\nu}}{f_{\nu,2}/z^{q_\nu}}
\]
in $0<\rho\le |z|$ in which the $f_{\nu,2}/z^{q_\nu}$ is clearly still analytic and not zero at the origin. We deduce from Lemmas \ref{L:lower-bound} and \ref{L:Miles} again that for $0\le\nu\le n-1$
	\begin{equation}\label{E:miles-2}
		\begin{split}
			\log| f_\nu| & =\log |f_{\nu,1}|-q_\nu\log r + \log |f_{\nu, 2}/z^{q_\nu}|^{-1}\\
                 & \le \log M(r,\, f_{\nu,1})-q_\nu\log r + H(\eta)\log M\big(2\,e\,r,\, f_{\nu,2}/z^{q_\nu}\big)\\
                   & \le \log^+ M(r,\, f_{\nu,1})-q_\nu\log r+ H(\eta)\log^+ M\big(2\,e\,r,\, f_{\nu,2}/z^{q_\nu}\big)\\
                 & \le \Big(\frac{R+r}{R-r}\Big)\,T(R,\,f_{\nu,\, 1})
                    -q_\nu\log r +H(\eta)\,\Big(\frac{R+2\,e\,r}{R-2\,e\,r}\Big)\, T(R,\, f_{\nu,2}/z^{q_\nu})\\ 
                 & \le  \Big(\frac{R+r}{R-r}\Big)\,T(R,\,f_{\nu,\, 1})
                    -q_\nu\log r\\
				& \qquad +H(\eta)\,\Big(\frac{R+2\,e\,r}{R-2\,e\,r}\Big)\, \big(T(R,\, f_{\nu,2})+q_\nu\log R\big)\\ 
                 & \le  B\Big(\frac{R+r}{R-r}\Big)\,T(CR,\,f_{\nu})
                    -q_\nu\log r \\
					&	\qquad +H(\eta)\,\Big(\frac{R+2\,e\,r}{R-2\,e\,r}\Big)\, \big(B T(CR,\, f_{\nu})+q_\nu\log R\big)\\ 
                 & <B\big(1+H(\eta)\big)\, \Big(\frac{R+2\,e\,r}{R-2\,e\,r}\Big)\,  T(CR,\, f_{\nu})+q_\nu\,\log\Big[\frac{R^{H(\eta)(\frac{R+2\,e\,r}{R-2\,e\,r})}}{r}\Big].
		\end{split}
	\end{equation}
 Applying (\ref{E:lower-bound}) to 
 (\ref{E:miles-1}) or (\ref{E:miles-2}) ($0\le\nu\le n-1$) and substituting them into (\ref{E:matrix}) completes the proof.
\qed
	
\bigskip

\section{Proof of Theorem \ref{T:density}}

We are now ready to prove Theorem \ref{T:density}. We choose $C=e$ in Theorem \ref{T:main}. Let $\alpha=2\, e$. We define annuli by
\[
\Lambda_j=\bigl\{ z: \alpha^j \le |z| \le \alpha^{j+3/2}\bigr\},\quad\quad
j=1,\, 2,\, \cdots.
\]
\medskip

We take $R=3\,er$ in (\ref{E:D}) in Theorem \ref{T:main} and suppose $z$ belongs to $\Omega\cap\Lambda_j$ where $\Omega$ is defined in Theorem \ref{T:main}. That is, we have, for each $0\le j\le n-1$, 
\begin{equation}
		\begin{split}
			m(r,\, y^{(j)}) & \le (2\pi+1)\,r\,n\Bigg\{1+(3e)^{5H(\eta)q}r^{(5H(\eta)-1)q}\,
\exp\Big[5B\big(1+H(\eta)\big)\, T(3\,e^2r)\Big]\Bigg\}\\
				&\quad +\log K_1\\
                 &= \log K_1+(2\pi+1)\, r\, D(f_\nu,\, \rho,\, r,\, 3er,\, \eta,\, B,\, e),
         \end{split}
	\end{equation}
where $D$ is given in (\ref{E:D}). Taking logarithm on both sides of the inequality once more yield the required estimate (\ref{E:density-bound}).
\bigskip

It remains to verify the size of the exceptional set of $r$, which follows from
\medskip

\begin{lemma}\sl
Let $\eta$ and $H(\eta)$ be defined in Lemma \ref{L:lower-bound} and Theorem \ref{T:main} respectively. Then the estimate (\ref{E:density-bound}) for a meromorphic solution $y(z)$ is valid for all $r$ sufficiently large except on a set of positive logarithmic density at most $16\eta e^{5/2}/(1+\log 2)$.
\end{lemma}
\bigskip

Let $E_j$ be the union of exceptional circles lying in $\Lambda_j$ and
\[
E(r)=[1,\, r)\cap\bigl(\cup_{j=1}^\infty E_j\bigr).
\]
Let $q=\left[\frac{\log r}{\log\alpha}\right]$, then Lemma \ref{L:lower-bound} gives
		\[
			\begin{split}
				\int_{E(r)}\frac{dt}{t} &\le\sum_{j=1}^q \int_{E_j}\frac{dt}{t}\\
                        & \le\sum_{j=1}^q \frac{4\eta(2\, er)}{\alpha^j}\\
                        & \le\sum_{j=1}^q \frac{4\eta(2e\alpha^{j+3/2})}{\alpha^j}\\
                        & \le\frac{\log r}{\log \alpha}(8\eta e\alpha^{3/2})\\
						& \le \Big(\frac{16\,\eta e^{5/2}}{1+\log 2}\Big)\,\log r.
			\end{split}
		\]
Thus
		\[
			\limsup_{r\to\infty} \frac{1}{\log r}\int_{E(r)}\frac{dt}{t}\le 
			\frac{16\eta e^{5/2}}{1+\log 2}<1.
		\]
\medskip

This completes the proof of the Lemma. This completes the proof of the theorem.\qed

\subsection*{Acknowledgment} The author would like to acknowledge Ms. Xudan Luo for her careful reading of the manuscript.
\medskip

\bibliographystyle{amsplain}

\end{document}